\documentclass[reqno]{amsart}
\usepackage{tikz,genyoungtabtikz,enumitem,rotating,latexsym,bm,stmaryrd}
 \usetikzlibrary{positioning,intersections}
  \usepackage{amsmath,amsthm,amsfonts,amssymb,mathrsfs,pb-diagram,hyperref}
 \usepackage{caption,cases}
\usepackage{lipsum,wasysym}
\usepackage{mathtools}
\usetikzlibrary{arrows,calc,through,backgrounds,matrix,decorations.pathmorphing}

\usepackage[vcentermath,enableskew]{youngtab}

\usepackage{cleveref}
\newtheorem{thm}{Theorem}[section]

\newtheorem*{prop*}{Proposition}
\newtheorem*{thm*}{Main Theorem}

\newtheorem*{cor*}{Corollary}
\newtheorem*{conj*}{Conjecture}

\newtheorem*{move1'}{Move $\mathbf{1^*}$}

\theoremstyle{remark}

\newtheorem{rmk}[thm]{Remark}

\newtheorem*{Acknowledgements*}{Acknowledgements}

\theoremstyle{definition}
\newtheorem{defn}[thm]{Definition}
\newtheorem{eg}[thm]{Example}

\crefname{defn}{Definition}{Definitions}
\crefname{thm}{Theorem}{Theorems}
\crefname{prop}{Proposition}{Propositions}
\crefname{lem}{Lemma}{Lemmas}
\crefname{cor}{Corollary}{Corollaries}
\crefname{conj}{Conjecture}{Conjectures}
\crefname{section}{Section}{Sections}
\crefname{subsection}{Subsection}{Subsections}
\crefname{eg}{Example}{Examples}
\crefname{figure}{Figure}{Figures}
\crefname{rem}{Remark}{Remarks}
\crefname{rmk}{Remark}{Remarks}
\crefname{equation}{equation}{equation}

\Crefname{defn}{Definition}{Definitions}
\Crefname{thm}{Theorem}{Theorems}
\Crefname{prop}{Proposition}{Propositions}
\Crefname{lem}{Lemma}{Lemmas}
\Crefname{cor}{Corollary}{Corollaries}
\Crefname{conj}{Conjecture}{Conjectures}
\Crefname{section}{Section}{Sections}
\Crefname{subsection}{Subsection}{Subsections}
\Crefname{eg}{Example}{Examples}
\Crefname{figure}{Figure}{Figures}
\Crefname{rem}{Remark}{Remarks}
\Crefname{rmk}{Remark}{Remarks}

\numberwithin{equation}{section}

\newcommand{\sts}{\mathsf{s}}  
\newcommand{\stt}{\mathsf{t}}  
\newcommand{\stu}{\mathsf{u}}  
\newcommand{\stv}{\mathsf{v}}  

 \newcommand{\SSTS}{\mathsf{S}}  
\newcommand{\SSTT}{\mathsf{T}}  
\newcommand{\SSTU}{\mathsf{U}} 
\newcommand{\SSTV}{\mathsf{V}}  

\newcommand{\down}{m{\downarrow}}
\newcommand{\up}{m{\uparrow}}
\newcommand{\Hom}{\operatorname{Hom}}

\newcommand{\suchthat}{\;\ifnum\currentgrouptype=16 \middle\fi|\;} 
\newcommand{\ZZ}{{\mathbb Z}}
\newcommand{\NN}{{{\mathbb Z}_{\geq0}}}

\newcommand{\CC}{\mathbb{Q}}

\newcommand\mptn[1]{\mathscr{P}_{#1}}

\usetikzlibrary{decorations.pathreplacing}

\usepackage{amsmath}
\mathchardef\mhyphen="2D

\DeclareMathOperator{\node}{node}

\newcommand{\Std}{\mathrm{Std}}

\newcommand{\SStd}{\mathrm{SStd}}
\newcommand{\Latt}{\mathrm{Latt}}

\let\originalleft\left
\let\originalright\right
\def\left#1{\mathopen{}\originalleft#1}
\def\right#1{\originalright#1\mathclose{}}

\renewcommand{\geq}{\geqslant}

\renewcommand{\leq}{\leqslant}

\def\ignore#1{\relax}

\def\ignore#1{\relax}

 \title[The lattice permutation condition  for Kronecker tableaux]{ The lattice permutation condition  \\ for Kronecker tableaux \\ (Extended abstract)}
\author{C. Bowman}
\author{M. De Visscher}
\author{J. Enyang}

\usepackage{scalefnt}
\begin{document}
 \maketitle
 \begin{abstract}
We recently generalised the lattice permutation condition for Young tableaux  to Kronecker tableaux and hence calculated a large new class of stable Kronecker coefficients labelled by co-Pieri triples.  
In this extended abstract we discuss   important families of co-Pieri triples for which our combinatorics 
simplifies drastically.  

\end{abstract}

\section{Introduction}

Perhaps the last major open problem in  the  complex representation theory of symmetric groups 
is to describe the decomposition of a tensor product of two simple representations.  
  The  coefficients describing the decomposition of these tensor products  are known 
   as the {\em Kronecker coefficients} and they have been described as 
    `perhaps the most challenging, deep and mysterious objects in algebraic combinatorics'.       Much recent progress has  focussed on   the stability properties enjoyed by  Kronecker coefficients.

\begin{figure}[ht!]$$  \begin{tikzpicture}[scale=0.4]
           \begin{scope}     \draw (0,3) node {  $\scalefont{0.5}\varnothing$  };   
  \draw (-3,0) node {   $ \frac{1}{2}$  };     \draw (-3,3) node {   $ 0$  };   
    \draw (-3,-3) node {   $ 1$  };     \draw (-3,-6) node {   $ 1\frac{1}{2}$  };   
    \draw (-3,-9) node {   $ 2$  };         \draw (-3,-12) node {   $ 2\frac{1}{2}$  };     
    \draw (-3,-15) node {   $ 3$  };     
              \draw (0,0) node { $\scalefont{0.5}\varnothing$  };   
    \draw (0,-3) node   {   \text{	$\scalefont{0.5}\varnothing$	}}		;    \draw (+3,-3) node   {  
     $ 
\scalefont{0.4}\yng(1)$	
    }		;
      \draw (0,-6) node   {   \text{	$\scalefont{0.5}\varnothing$	}};
     \draw (3,-6) node   {  $       \,
\scalefont{0.5}\yng(1) $	 };
    \draw[<-] (0.0,1) -- (0,2);     \draw[->] (0.0,-0.75) -- (0,-2.25);        
                \draw[->] (0.5,-0.75) -- (2.5,-2.25);    \draw[->] (2.5,-3.75) -- (0.5,-5.25);       \draw[->] (0,-3.75) -- (0,-5.25);   \draw[->] (3,-3.75) -- (3,-5.25);
  \draw[->] (0,-6.75) -- (0,-8.25);   \draw[->] (03,-6.75) -- (3,-8.25); 
    \draw[->] (0.5,-6.75) -- (2.5,-8.25);   \draw[->] (4,-6.75) -- (8,-8.25);   \draw[->] (3.5,-6.75) -- (5,-8.25); 
     \draw (+0,-9) node   {   	  $\scalefont{0.5}\varnothing$  };                 
          \draw (+3,-9) node   {  $ 
\scalefont{0.4}\yng(1) $	 }		;
             \draw (+6,-9) node              {    $
\scalefont{0.4}\yng(2)$	 	 }		;
                          \draw (+9,-9) node            {    $
\scalefont{0.4}\yng(1,1)$	 	 }			;
   \draw[<-] (0,-11.25) -- (0,-9.75);   \draw[<-] (03,-11.25) -- (3,-9.75);  \draw[<-] (06,-11.25) -- (6,-9.75); \draw[<-] (9,-11.25) -- (9,-9.75); 
    \draw[<-] (0.5,-11.25) -- (2.5,-9.75);   \draw[<-] (4,-11.25) -- (8,-9.75);   \draw[<-] (3.5,-11.25) -- (5,-9.75); 
     \draw (+0,-12) node   {   	$\scalefont{0.5}\varnothing$  };                 
          \draw (+3,-12) node   {  $ 
\scalefont{0.4}\yng(1)$	 }		;
             \draw (+6,-12) node              {    $
\scalefont{0.4}\yng(2)$	 	 }		;
                          \draw (+9,-12) node            {    $
\scalefont{0.4}\yng(1,1)$	 	 }			;
   \draw[->] (0,-12.75) -- (0,-14.25);   \draw[->] (03,-12.75) -- (3,-14.25);  \draw[->] (06,-12.75) -- (6,-14.25); \draw[->] (9,-12.75) -- (9,-14.25); 
    \draw[->] (0.5,-12.75) -- (2.5,-14.25);   \draw[->] (4,-12.75) -- (8,-14.25);   \draw[->] (3.5,-12.75) -- (5,-14.25); 
     \draw (+0,-15) node   {   	$\scalefont{0.5}\varnothing$  };                 
          \draw (+3,-15) node   {  $ 
\scalefont{0.4}\yng(1) $	 }		;
             \draw (+6,-15) node              {    $
\scalefont{0.4}\yng(2)$	 	 }		;
                          \draw (+9,-15) node            {    $
\scalefont{0.4}\yng(1,1)$	 	 }			;
  \draw[->] (6.5,-12.75) -- (12,-14.25);    \draw[->] (7,-12.75) -- (14.75,-14.25);  
    \draw[->] (9.5,-12.75) -- (15.25,-14.25);    \draw[->] (10,-12.75) -- (17.7,-14.1);  
     \draw (12,-15) node            {    $
\scalefont{0.4}\yng(3)$	 	 }			;
       \draw (15,-15) node            {    $
\scalefont{0.4}\yng(2,1)$	 	 }			;
       \draw (18,-15) node            {    $
\scalefont{0.4}\yng(1,1,1)$	 	 }			;
    \end{scope}\end{tikzpicture} 
     $$
  
     \caption{The first three layers of the branching graph $\mathcal{Y}$}
     \label{oscillate}\label{brancher}
\end{figure}
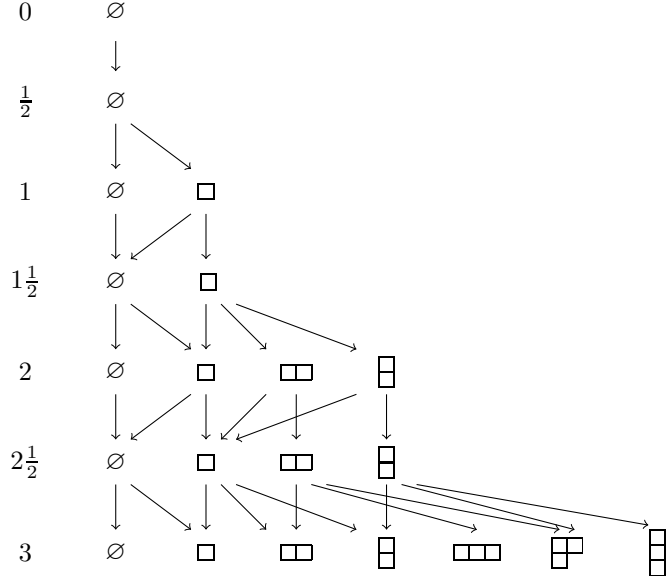

  Whilst a complete understanding of the  Kronecker coefficients seems out of reach, 
the purpose of this work is to attempt to understand
 the  {\em stable} Kronecker coefficients in terms of oscillating tableaux.
 Oscillating tableaux  hold a distinguished position in the study of 
  tensor product decompositions   \cite{MR1035496,MR3090983,MR2264927}  but surprisingly they  have never before been 
  used to calculate Kronecker coefficients of symmetric groups.  
 In this work, we  see that the oscillating tableaux defined as paths on the graph given in   \cref{brancher} (which we call Kronecker tableaux) provide
  bases of certain modules for the partition algebra, $P_s(n)$, which is closely related to the symmetric group.
 We hence add a new level of structure to the classical picture --- this extra structure is the key to our main result: the co-Pieri rule for stable Kronecker coefficients.

A momentary glance at the   graph given in \cref{oscillate}
  reveals a very familiar subgraph: namely Young's graph (with each level doubled up).   
 The stable Kronecker coefficients labelled by triples from  this subgraph are  well-understood  --- 
  the values of these coefficients
    can be calculated via a tableaux counting algorithm known as the Littlewood--Richardson rule \cite{LR34}.   
  This rule    has long served as the hallmark for our understanding  of Kronecker coefficients.  
The Littlewood--Richardson  rule was discovered as a rule of two halves (as we explain below).  
In \cite{BDE} we succeed in  generalising one half of  this rule to all Kronecker tableaux,  and thus solve one half of the stable Kronecker problem.  
 Our main result unifies and vastly generalises the work of Littlewood--Richardson \cite{LR34} and many other authors \cite{RW94,Rosas01,ROSAANDCO,BWZ10,MR2550164}.  
 Most promisingly, our   result counts  explicit   homomorphisms 
and thus  works on a structural level above any description of a family of Kronecker coefficients   since those first considered by Littlewood--Richardson   \cite{LR34}.

In more detail,  given a triple of partitions $(\lambda,\nu, \mu)$  and with $|\mu|=s$, we have an associated 
  skew $P_s(n)$-module spanned by the Kronecker  tableaux  from $\lambda$ to $\nu$  of length $s$, which
  we denote by $\Delta_s(\nu \setminus\lambda )   $.  
  For $\lambda=\scalefont{0.5}\varnothing$ and $n\geq 2s$  these 
  modules provide a complete set of non-isomorphic $P_s(n)$-modules (and we drop the 
  partition $\scalefont{0.5}\varnothing$ from the notation).    
  The stable Kronecker coefficients are then interpreted as the dimensions, 
\begin{equation}\label{dagger}\tag{$\dagger$}
\overline{g}(\lambda,\nu,\mu)
=  \dim_\CC( \Hom_{  P_{s}(n)}( \Delta_{s}(\mu), \Delta_s(\nu \setminus\lambda )    ) )    
\end{equation}for $n\geq 2s$.  
Restricting to the Young subgraph, or equivalently  to a triple     $(\lambda,\nu,\mu)$  of so-called {\em maximal depth} such that   $|\lambda| + |\mu| = |\nu|$, 
 these modules specialise to   the usual simple and skew modules for symmetric groups;  hence the multiplicities $\overline{g}(\lambda,\nu,\mu)$ are the Littlewood--Richardson coefficients. We hence recover the well-known fact that the Littlewood--Richardson coefficients appear as the subfamily of stable Kronecker coefficients labelled by triples of maximal depth.
 The tableaux   counted by the Littlewood--Richardson  rule satisfy 2 conditions: the {\em semistandard} and  {\em lattice permutation} conditions. 
%
%
%
%
 In 
\cite{BDE} we generalise the lattice permutation condition to 
  Kronecker tableaux. 

 \begin{thm*}[{\cite[Main Theorem]{BDE}}]
Let $(\lambda,\nu,\mu)$ be a    {co-Pieri} triple     or a triple of maximal depth. Then the stable Kronecker coefficient $\overline{g}(\lambda, \nu, \mu)$ is given by the number of  
 semistandard  Kronecker tableaux of shape $\nu\setminus\lambda$ and weight $\mu$ whose reverse reading word is a lattice permutation.   
 \end{thm*}
The observant reader will notice that the statement above describes the Littlewood--Richardson coefficients uniformly as part of a far broader family of stable Kronecker coefficients (and is the first result in the literature to do so). Whilst the classical Pieri rule (describing the semistandardness condition for Littlewood--Richardson tableaux) is elementary, it served as a first step towards 
   understanding the full Littlewood--Richardson rule; indeed 
   Knutson--Tao--Woodward  have shown that  the  Littlewood--Richardson rule follows from the Pieri rule by associativity  \cite{taoandco}.  
    We  hope that our  generalisation of the   co-Pieri rule (the lattice permutation condition for Kronecker tableaux)  will prove equally useful in the study of  stable   Kronecker coefficients. 

   The definition of {\sf semistandard Kronecker tableaux} naturally generalises the classical notion of semistandard Young tableaux as certain ``orbits" of paths on the branching graph given in \cref{brancher} (see Section 1.2 and Definition \ref{sstrd}). The {\sf lattice permutation condition} is identical to the classical case once we generalise the dominance order to all steps in the branching graph $\mathcal{Y}$ to define the reverse reading word of a semistandard Kronecker tableau  (see \cref{sec:latticed}). 
   
\bigskip

\noindent \textbf{Examples of co-Pieri triples.} The definition of {\em co-Pieri triples } is given in \cite[Theorem 4.12]{BDE} and can appear quite technical at first reading; we present a few special cases here.   
\begin{itemize}[leftmargin=*]
\item[$(i)$] $\lambda$ and $\mu$ are one-row partitions and $\mu$ is arbitrary.
This family has been extensively studied over the past thirty years and 
there are   many 
 distinct  combinatorial descriptions of some or all of these coefficients 
 \cite{RW94,Rosas01,ROSAANDCO,BWZ10,MR2550164}, none of which generalises.

\item[$(ii)$] the two skew partitions $\lambda \ominus (\lambda \cap \nu)$ and $\nu \ominus (\lambda \cap \nu)$ have no two boxes in the same column and   $|\mu| = \max \{|\lambda \ominus (\lambda \cap \nu)| , |\nu \ominus (\lambda \cap \nu)|\}$. 
 It is easy to see that if, in addition, $(\lambda, \nu, \mu)$ is a triple of maximal depth, then this case specialises to the classical co-Pieri triples.

\item[$(iii)$] $\lambda = \nu =  (dl,d(l-1), \ldots , 2d,d)$ for any $l,d\geq 1$ and $|\mu| \leq d$.
\end{itemize}

In this extended abstract we have chosen to focus primarily on case $(i)$ as these triples  carry many of the tropes of  general co-Pieri triples (but with significant simplifications which serve to make this abstract more approachable) and because case $(i)$ should be familiar to many readers due to its many appearances in the literature.

 \section{The partition algebra and Kronecker tableaux }\label{sec2} 
   \label{sec:standard}

    The combinatorics underlying the representation theory of the partition algebras and   symmetric 
     groups  is based on     partitions.  
A {\em  partition}    $\lambda $ of $n$, denoted $\lambda \vdash n$, is defined to be a  sequence   of weakly decreasing non-negative integers which  sum to $n$.  
       We let $\varnothing$  denote the unique partition of 0.  
 Given a partition, $\lambda=(\lambda _1,\lambda _2,\dots )$, the  associated   {\em Young diagram}  is the set of nodes
$[\lambda]=\left\{(i,j)\in\mathbb{Z}_{>0}^2\ \left|\ j\leq \lambda_i\right.\right\}.$ 
We define the length, $\ell(\lambda)$, of a partition $\lambda$, to be the number of non-zero parts.   
Given  $\lambda = (\lambda_1,\lambda_2, \ldots,\lambda_{\ell} )$  a partition and $n$   an integer, define 
   $\lambda_{[n]}=(n-|\lambda|, \lambda_1,\lambda_2, \ldots,\lambda_{\ell}).$   
 Given $\lambda_{[n]} $ a partition of $n$, we say that the partition has {\em depth} equal to $|\lambda|$.

    The partition algebra  is generated  as an algebra  by the elements $s_{k,k+1}$, $p_{k+1/2}$ ($1\leq k\leq r-1$) and $p_k$ ($1\leq k\leq r$) pictured below modulo a long list of relations.  One can visualise any product in this algebra as simply being given by concatenation of diagrams, modulo some surgery to remove closed loops \cite{BDE}.  
$$
{s}_{k,k+1}=
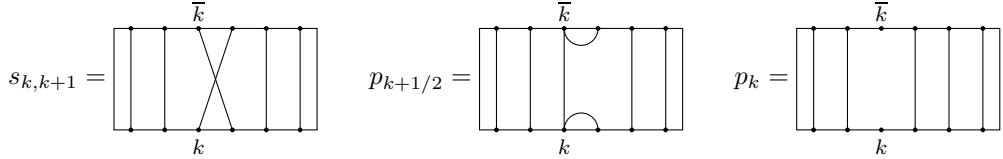
\begin{minipage}{34mm}\scalefont{0.8}\begin{tikzpicture}[scale=0.45]
  \draw (0,0) rectangle (6,3);
  \foreach \x in {0.5,1.5,...,5.5}
    {\fill (\x,3) circle (2pt);
     \fill (\x,0) circle (2pt);}
     \draw (2.5,-0.49) node {$k$};
    \draw (2.5,+3.5) node {$\overline{k}$};
  \begin{scope}
    \draw (0.5,3) -- (0.5,0);
    \draw (5.5,3) -- (5.5,0);
        \draw (2.5,3) -- (3.5,0);
            \draw (4.5,3) -- (4.5,0);
    \draw (3.5,3) -- (2.5,0);
    \draw (1.5,3) -- (1.5,0);
   \end{scope}
\end{tikzpicture}\end{minipage}    
p_{k+1/2}
=\begin{minipage}{34mm}\scalefont{0.8}\begin{tikzpicture}[scale=0.45]
  \draw (0,0) rectangle (6,3);
  \foreach \x in {0.5,1.5,...,5.5}
    {\fill (\x,3) circle (2pt);
     \fill (\x,0) circle (2pt);}
     \draw (2.5,-0.49) node {$k$};
    \draw (2.5,+3.5) node {$\overline{k}$};
  \begin{scope}
    \draw (0.5,3) -- (0.5,0);
    \draw (5.5,3) -- (5.5,0);
     \draw (1.5,3) -- (1.5,0);
\draw (3.5,0) arc (0:180:.5 and 0.5);
\draw (2.5,3) arc (180:360:.5 and 0.5);
    \draw (2.5,3) -- (2.5,0);
    \draw (4.5,3) -- (4.5,0);
   \end{scope}
\end{tikzpicture}\end{minipage}  
p_k 
=\begin{minipage}{34mm}\scalefont{0.8}\begin{tikzpicture}[scale=0.45]
  \draw (0,0) rectangle (6,3);
  \foreach \x in {0.5,1.5,...,5.5}
    {\fill (\x,3) circle (2pt);
     \fill (\x,0) circle (2pt);}
   
    \draw (2.5,-0.49) node {$k$};
    \draw (2.5,+3.5) node {$\overline{k}$};
  \begin{scope}
    \draw (0.5,3) -- (0.5,0);
    \draw (5.5,3) -- (5.5,0);
     \draw (4.5,0) -- (4.5,3); 
       
    \draw (1.5,3) -- (1.5,0);
    \draw (3.5,3) -- (3.5,0);
   \end{scope}
\end{tikzpicture}\end{minipage}  
 $$

   Define    the branching graph $\mathcal{Y}$ as follows.   For $k\in \NN$, we denote by $\mathscr{P}_{\leq k}$ the set of partitions of degree less or equal to $k$. Now the set of vertices on the $k$th and $(k+1/2)$th levels of $\mathcal{Y}$ are given by 
    $$  
   {\mathcal{Y}}_{k}=  \{ (\lambda,k-|\lambda|) \mid \lambda \in \mathscr{P}_{\leq k}\}
   \qquad
   {\mathcal{Y}}_{k+1/2} = 
\{ (\lambda,k-|\lambda|) \mid \lambda \in \mathscr{P}_{\leq k}\}.$$   The edges of $\mathcal{Y}$ are   as follows,
\begin{itemize}[leftmargin=*] 
\item for $(\lambda,l) \in {\mathcal{Y}}_k$ and $(\mu,m) \in {\mathcal{Y}}_{k+1/2}$ there is an edge $(\lambda,l) \to(\mu,m)$ 
 if $\mu = \lambda$,  or   
if $\mu $ is  obtained from $\lambda $ by removing a box in the $i$th row for some $i\geq 1$;
  we write $\mu =\lambda- \varepsilon_0$ or $\mu =\lambda- \varepsilon_i$, respectively. 
\item for $(\lambda,l)  \in  {\mathcal{Y}}_{k+1/2}$ and $(\mu,m)  \in {\mathcal{Y}}_{k+1}$ there is an edge  $(\lambda,l) \to(\mu,m)$ 
if $\mu = \lambda$,  or   
if    $\mu $ is   obtained from $\lambda $ by adding a box in the $i$th row for some $i\geq 1$;
  we write $\mu =\lambda+ \varepsilon_0$ or $\mu =\lambda+ \varepsilon_i$, respectively. 
 \end{itemize}
When it is convenient, we  decorate each edge with the index of the node that is added or removed when reading down the diagram.  
The first few levels
of $\mathcal{Y}$  are given in   Figure \ref{brancher}.  
 When no confusion is possible, we identify $(\lambda,l) \in \mathcal{Y}_{k}$ with the partition $\lambda$.  
 \begin{defn}
Given   $\lambda  \in \mathscr{P}_{r-s} \subseteq \mathcal{Y}_{r-s}$ and $\nu \in \mathscr{P}_ {\leq r} \subseteq \mathcal{Y}_{r}$, we  
define a  {\em  standard Kronecker tableau} of shape $ \nu \setminus  \lambda $ and degree $s$  to be a path $\stt$  of the form 
\begin{equation}\label{genericpath} 
\lambda = \stt(0) \to \stt(\tfrac{1}{2}) \to    \stt(1)\to  \dots \to  \stt(s-\tfrac{1}{2})\to \stt(s) = \nu, 
 \end{equation}
 in other words $\stt$ is  a path in $\mathcal{Y}$ which begins at $\lambda$ and terminates at $\nu$.  
We let $\Std_s(\nu \setminus  \lambda)$ denote the set of all such paths.  If $\lambda = \emptyset \in \mathcal{Y}_0$ then we write $\Std_r(\nu)$ instead of $\Std_r(\nu \setminus \emptyset)$.  
Given $\sts,\stt $ two standard Kronecker tableaux of degree $s$, we write $\sts \trianglerighteq \stt $ if 
$\sts(k)\trianglerighteq \stt(k)$ for all $0\leq k\leq s$.
\end{defn}

We can think of a path as either the sequence of partitions or the sequence of boxes removed and added.  
We usually prefer the latter case and record these boxes removed and added pairwise. For a pair
 $(-\varepsilon_p,+\varepsilon_q)$ we call this an add or remove step if $p=0$ or $q=0$ respectively 
 (because the effect of this step is to add or remove a box) and we call this a dummy step if $p=q$ (as we end up at the same partition as we started);  
 we write $a(q)$ or $r(p)$ for an add or remove step and $d(p)$ for a dummy step.  
Many examples are given below, in particular the reader should compare the paths of Example \ref{example3} with those depicted in the central diagram in Figure \ref{maximaldepth}.   
 We let $\stt^\lambda$ denote the most dominant element of $\Std_s(\lambda)$, namely that of the form:
$$
\underbrace{
d(0) \circ d(0) \circ 
\dots
\circ d(0) }_{r-|\lambda|}
\circ 
\underbrace{
a(1)
\circ 
\dots\circ 
a(1)}_{\lambda_1}\circ 
\underbrace{
a(2)
\circ 
\dots\circ 
a(2)}_{\lambda_2}\circ 
\cdots
 $$ 
 Given $\lambda
\in\mathscr{P}_{r-s} \subseteq \mathcal{Y}_{r-s}$ and  $\nu \in\mathscr{P}_{\leq r} \subseteq \mathcal{Y}_r$, define
  the {\em skew cell module} 
 $$
\Delta_s(\nu\setminus\lambda) = {\rm Span}\{ \stt^\lambda \circ \sts \mid \sts \in \Std_s(\nu\setminus\lambda)\}
$$
with the action of $P_s(n)\hookrightarrow P_{r-s}(n) \otimes P_s(n)\hookrightarrow P_r(n)$  given as in \cite[Section 2.3]{BDE}. 
 If $\lambda =\varnothing$, then we simply denote this module by $\Delta_s(\nu)$.  
      Let $\lambda\in \mathscr{P}_{r-s}$, $\mu\in \mathscr{P}_s$ and $\nu \in \mathscr{P}_{\leq r}$. Then 
  we are able to define the stable Kronecker coefficients (even if this is not their usual definition)  to be the multiplicities
$$
\overline{g}(\lambda,\nu,\mu)
 =  \dim_\CC( \Hom_{  P_{s}(n)}( \Delta_{s}(\mu), \Delta_s(\nu \setminus\lambda )    ) ) 
  $$
  for all $n\geq 2s$.   When $s=|\nu|-|\lambda|$, the (skew) cell modules for partition algebras specialise to the usual Specht modules of the symmetric groups and 
we hence easily see that these stable coefficients coincide with the classical Littlewood--Richardson coefficients.

\section{The action of the partition algebra}  
  
 Understanding the action of the partition algebra on skew modules is difficult in general.  
 In this section, we show that this can be done to some extent in the cases of interest to us.  
  We have assumed that $|\mu|=s$, therefore  
  the ideal $P_s(n)p_r P_s(n)\subset P_s(n)$  annihilates $\Delta_s(\mu)$ and this  motivates the following definition.

 \begin{defn} We define the {\sf Dvir radical} of the  skew module $\Delta_s(\nu\setminus\lambda)$ by 
$$  {\sf DR}_s(\nu\setminus\lambda)  =  \Delta_s(\nu\setminus\lambda)P_s(n)p_rP_s(n) 
\subseteq \Delta_s(\nu\setminus\lambda) $$
and set $$\Delta^0_s(\nu\setminus\lambda)= 
 \Delta_s(\nu\setminus\lambda) /{\sf DR}_s(\nu\setminus\lambda). $$
  If $s=|\nu|-|\lambda|$, then    set 
  $  \Std^0_s(\nu\setminus\lambda) =   \Std_s(\nu\setminus\lambda).
  $
  If $\lambda$ and $\nu$ are one-row partitions, then    set 
 $   \Std^0_s(\nu\setminus\lambda) \subseteq   \Std_s(\nu\setminus\lambda) 
 $ 
  to be the  subset of paths, $\sts$, whose steps are of the form
  $$
  r(1)=(-1,+0) \qquad d(1)=(-1,+1)  \qquad a(1)=(-0,+1)
  $$
  and such that the total number of boxes removed in $\sts$ is less than or equal to $|\lambda|$.  
  \end{defn}


 Fix  $\stt \in   {\Std}_r( \nu   )$ and $1\leq k \leq r$  and suppose that
    $$ \stt{(k-1)}    \xrightarrow{-t}  \stt(k-\tfrac{1}{2})    \xrightarrow{+u}   \stt(k+1) \xrightarrow{-v}      \stt(k+\tfrac{1}{2})  \xrightarrow{+w}  \stt(k+1).$$
We define  $ \stt_{k \leftrightarrow k+1}\in \Std_r(\nu)$ to be the tableau, if it exists, determined by  $  \stt_{k \leftrightarrow k+1}(l) =\stt(l) $ for $l\neq k, k \pm \tfrac{1}{2}  $ and 
 $$  \stt_{k \leftrightarrow k+1} {(k-1)}    \xrightarrow{-v}     \stt_{k \leftrightarrow k+1}{(k-\tfrac{1}{2})}   \xrightarrow{+w}  
   \stt_{k \leftrightarrow k+1}{(k)}  \xrightarrow{-t}       \stt_{k \leftrightarrow k+1}(k+\tfrac{1}{2})  \xrightarrow{+u}  \stt_{k \leftrightarrow k+1}(k+1).$$
 Let   $(\lambda,\nu,s)  $  be  such that 
 $s=|\nu|-|\lambda|$, or $\lambda$ and $\nu$ are both one-row partitions, then $\Delta^0_s(\nu\setminus\lambda)$ is free as a $\ZZ$-module with  basis 
$$  \{\stt \mid \stt \in \Std^0_s(\nu\setminus\lambda)\} $$
and  the $P_s(n)$-action on $\Delta_s^0(\nu\setminus\lambda)$ is   as follows:
\begin{equation}\label{co-case}
(\stt + {\sf DR_s(\nu\setminus\lambda)} ) s_{k,k+1}  =
\begin{cases}
  {\stt_{k\leftrightarrow k+1}}  + {\sf DR_s(\nu\setminus\lambda)} &\text{if  ${\stt_{k\leftrightarrow k+1}}$ exists}	\\
  - \stt +\sum_{\sts\rhd \stt} r_{\sts\stt}\sts+ {\sf DR_s(\nu\setminus\lambda)} &\text{otherwise}
  \end{cases}
\end{equation}
for $1\leq k < s$ and 
 $  (\stt + {\sf DR_s(\nu\setminus\lambda)} ) p_{k,k+1}  =
0$   and $
 ( \stt + {\sf DR_s(\nu\setminus\lambda)} ) p_{k}  =
0
 $
for   $1\leq k   \leq s$.  
The coefficients  $r_{\sts\stt}\in \mathbb{Q}$ are given in \cite[Theorem 2.9]{BDE}.

\begin{eg}
The set $\Std^0((3,3)\setminus (2,1))$ has two elements
$$
\stt_1=a(1) \circ a(2)\circ a(2)
\qquad
\stt_2=a(2) \circ a(1)\circ a(2).
$$
These are depicted on the lefthand-side of \ref{actioner}.  
We have that 
$$
s_{1,2}=\left(\begin{array}{cc}0 & 1 \\1 & 0\end{array}\right)
\qquad
s_{2,3}=\left(\begin{array}{cc}1 & -1 \\0 & -1\end{array}\right).$$
   \end{eg}  
  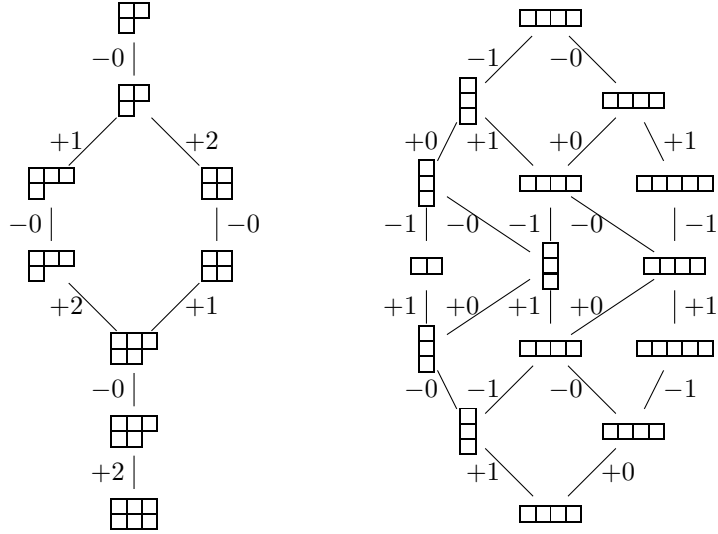
\begin{figure}[ht!]  $$  
       \begin{tikzpicture}[scale=0.55]
  
    \path   (0,0)edge[decorate]  node[left] {$-0$}  (0,-2);
          \path   (0,-2)edge[decorate]  node[right] {$+2$}  (2,-4);
                    \path   (0,-2)edge[decorate]  node[left] {$+1$}  (-2,-4);

                    \path   (-2,-6)edge[decorate]  node[left] {$-0$}  (-2,-4);                         
          \path   (2,-6)edge[decorate]  node[right] {$-0$}  (2,-4);
                                       
                  \path   (-2,-6)edge[decorate]  node[left] {$+2$}  (0,-8);                         
          \path   (2,-6)edge[decorate]  node[right] {$+1$}  (0,-8);

       \path   (0,0-8)edge[decorate]  node[left] {$-0$}  (0,-2-8);
                     \path   (0,-2-8)edge[decorate]  node[left] {$+2$}  (0,-4-8);
                         
%
%

%
%
       
%
%
                  \fill[white] (0,0) circle (17pt);   
  \begin{scope}   
 \fill[white] (0,0) circle (17pt);
     \draw (0,0) node {$  \scalefont{0.4}\yng(2,1) $  };   
\fill[white] (0,-2) circle (17pt);    \draw (0,-2) node{$  \scalefont{0.4}\yng(2,1)$  };
\fill[white] (-2,-4) circle (17pt);  \draw (-2,-4) node{$  \scalefont{0.4}\yng(3,1) $  };
\fill[white] (2,-4) circle (17pt);  \draw (2,-4) node{$   \scalefont{0.4}\yng(2,2) $  };
  \fill[white] (-2,-6) circle (17pt);  \draw (-2,-6) node{$  \scalefont{0.4}\yng(3,1)$  };
\fill[white] (2,-6) circle (17pt);  \draw (2,-6) node{$  \scalefont{0.4}\yng(2,2) $  };
  \fill[white] (0,-8) circle (17pt);  \draw (0,-8) node{$  \scalefont{0.4}\yng(3,2)$  };
 \fill[white] (0,-2-8) circle (17pt);    \draw (0,-2-8) node{$  \scalefont{0.4}\yng(3,2)$  };
\fill[white] (0,-4-8) circle (17pt);  \draw (0,-4-8) node{$  \scalefont{0.4}\yng(3,3) $  };

      \end{scope}
   \end{tikzpicture} \qquad\qquad
     \begin{tikzpicture}[scale=0.55]
  
  \path   (0,0)edge[decorate]  node[left] {$-1$}  (-2,-2);
    \path   (0,0)edge[decorate]  node[left] {$-0$}  (2,-2);
    
  \path   (-3,-4)edge[decorate]  node[left] {$+0$}  (-2,-2);    
  \path   (0,-4)edge[decorate]  node[left] {$+1$}  (-2,-2);    
  
    \path   (3,-4)edge[decorate]  node[right] {$+1$}  (2,-2);    
  \path   (0,-4)edge[decorate]  node[left] {$+0$}  (2,-2);

  \path   (-3,-8)edge[decorate]  node[left] {$-0$}  (-2,-10);    
  \path   (0,-8)edge[decorate]  node[left] {$-1$}  (-2,-10);

  \path   (3,-8)edge[decorate]  node[right] {$-1$}  (2,-10);    
  \path   (0,-8)edge[decorate]  node[left] {$-0$}  (2,-10);

  \path   (-3,-4)edge[decorate]  node[left] {$-1$}  (-3,-6);    
  \path   (-3,-4)edge[decorate]  node[left] {$-0$}  (0,-6);

  \path   (0,-4)edge[decorate]  node[left] {$-0$}  (3,-6);    
  \path   (0,-4)edge[decorate]  node[left] {$-1$}  (0,-6);      
  
  \path   (3,-4)edge[decorate]  node[right] {$-1$}  (3,-6);

  \path   (-3,-8)edge[decorate]  node[left] {$+1$}  (-3,-6);    
  \path   (-3,-8)edge[decorate]  node[left] {$+0$}  (0,-6);

  \path   (0,-8)edge[decorate]  node[left] {$+0$}  (3,-6);    
  \path   (0,-8)edge[decorate]  node[left] {$+1$}  (0,-6);      
  
  \path   (3,-8)edge[decorate]  node[right] {$+1$}  (3,-6);     
            \path   (2,-2-8)edge[decorate]  node[right] {$+0$}  (0,-4-8);
                    \path   (-2,-2-8)edge[decorate]  node[left] {$+1$}  (0,-4-8);
 
%
%
%
%
%
%
%
%
%
%
%
                  \fill[white] (0,0) circle (17pt);   
  \begin{scope}   
 \fill[white] (0,0) circle (18pt);
     \draw (0,0) node {$  \scalefont{0.4}\yng(4) $  };   
\fill[white] (-2,-2) circle (17pt);    \draw (-2,-2) node{$  \scalefont{0.4}\yng(1,1,1)$  };
\fill[white] (2,-2) circle (17pt);    \draw (2,-2) node{$  \scalefont{0.4}\yng(4)$  };

  \fill[white] (-3,-4) circle (17pt);  \draw (-3,-4) node{$  \scalefont{0.4}\yng(1,1,1)$  };
  \fill[white] (0,-4) circle (17pt);  \draw (0,-4) node{$  \scalefont{0.4}\yng(4)$  };
  \fill[white] (3,-4) circle (17pt);  \draw (3,-4) node{$  \scalefont{0.4}\yng(5)$  };

  \fill[white] (-3,-6) circle (17pt);  \draw (-3,-6) node{$  \scalefont{0.4}\yng(2)$  };
  \fill[white] (0,-6) circle (17pt);  \draw (0,-6) node{$  \scalefont{0.4}\yng(1,1,1)$  };
  \fill[white] (3,-6) circle (17pt);  \draw (3,-6) node{$  \scalefont{0.4}\yng(4)$  };

  \fill[white] (-3,-8) circle (17pt);  \draw (-3,-8) node{$  \scalefont{0.4}\yng(1,1,1)$  };
  \fill[white] (0,-8) circle (17pt);  \draw (0,-8) node{$  \scalefont{0.4}\yng(4)$  };
  \fill[white] (3,-8) circle (17pt);  \draw (3,-8) node{$  \scalefont{0.4}\yng(5)$  };

%
%
%
%
%
%
%
%
%
%
\fill[white] (-2,-2-8) circle (17pt);    \draw (-2,-2-8) node{$  \scalefont{0.4}\yng(1,1,1)$  };
\fill[white] (2,-2-8) circle (17pt);    \draw (2,-2-8) node{$  \scalefont{0.4}\yng(4)$  };
 \fill[white] (0,-4-8) circle (17pt);  \draw (0,-4-8) node{$   \scalefont{0.4}\yng(4) $  };
%
%
%

      \end{scope}
   \end{tikzpicture} 
    $$
    \caption{Oscillating  tableaux of shape $(3,3)\setminus (2,1)$ and $(4)\setminus (4)$ and degree 3. }
    \label{actioner}
    \end{figure}

        \begin{eg}\label{example3}
The set $\Std_3^0((4)\setminus(4))$ consists of the 7 oscillating tableaux
$$\begin{array}{ccccc}
\sts_1= r(1)\circ d(1)\circ a(1)  &
\sts_2= d(1)\circ r(1)\circ a(1)  &
\sts_3=  r(1)\circ a(1)  \circ d(1)   
\\
\sts_4= a(1)\circ r(1)\circ d(1)  &
\sts_5=  d(1)\circ a(1)  \circ r(1)   &
\sts_6=  a(1)\circ d(1)\circ r(1)   &
\\
    &   \sts_7=d(1)\circ d(1)  \circ d(1)
       \end{array}$$
  pictured in \cref{actioner}.     We have that 
$$
s_{1,2}=
\left(\begin{array}{ccccccc}
\cdot  & 1 	&	 \cdot  & \cdot  & \cdot  & \cdot  & \cdot  \\
1	 & \cdot  & \cdot  & \cdot  & \cdot  & \cdot  & \cdot  \\
\cdot  & \cdot  & \cdot  & 1  & \cdot  & \cdot  & \cdot  \\
\cdot  & \cdot  & 1  & \cdot  & \cdot  & \cdot  & \cdot  \\
\cdot  & \cdot  & \cdot  & \cdot  & \cdot  & 1  & \cdot  \\
\cdot  & \cdot  & \cdot  & \cdot  & 1  & \cdot  & \cdot  \\
\cdot  & \cdot  & \cdot  & \cdot  & \cdot  & \cdot  & 1  \\\end{array}\right)
\qquad
s_{2,3}=
\left(\begin{array}{ccccccc}
\cdot  & \cdot  & 1  	& \cdot  & \cdot  & \cdot  & \cdot  \\
\cdot  & \cdot  & \cdot  & \cdot  & 1 & \cdot  & \cdot  \\
1	  & \cdot  & \cdot  & \cdot  & \cdot  & \cdot  & \cdot  \\
\cdot  & \cdot  & \cdot  & \cdot  & \cdot  & 1 & \cdot  \\
\cdot  & 1  & \cdot  & \cdot  & \cdot  & \cdot  & \cdot  \\
\cdot  & \cdot  & \cdot  & 1  & \cdot  & \cdot  & \cdot  \\
\cdot  & \cdot  & \cdot  & \cdot  & \cdot  & \cdot  & 1 \\\end{array}\right)
$$
It is not difficult to see that this module decomposes as follows
$$
\Delta^0_3((4)\setminus(4))= 2 \Delta^0_3((3))
\oplus  2 \Delta^0_3((2,1))
\oplus \Delta^0_3((1^3)).
$$
\end{eg}

\section{Semistandard Kronecker tableaux}\label{sec:semistandard}

For any $(\lambda,\nu,s) \in \mptn{r-s}\times  \mptn{\leq r} \times  \ZZ_{> 0} $ and any $\mu \vdash s$ we have
 $$
\overline g( \lambda, \nu,\mu) = \dim_\CC \Hom_{P_s(n)}(\Delta_s(\mu), \Delta_s^0(\nu\setminus\lambda) ) = \dim_{\CC} \Hom_{\CC \mathfrak{S}_s}({\sf S}(\mu), \Delta_s^0(\nu \setminus \lambda)),
$$
where $\CC \mathfrak{S}_s$ is viewed as the quotient of $P_s(n)$ by the ideal generated by $p_r$.  
Now for each   $\mu = (\mu_1, \mu_2, \ldots , \mu_l) \vdash s$ we have an associated   Young permutation module 
 ${\sf M}  (\mu) = \CC\otimes_{\mathfrak{S}_\mu} \CC\mathfrak{S}_s$ 
 where $\mathfrak{S}_\mu = \mathfrak{S}_{\mu_1} \times \mathfrak{S}_{\mu_2}\times \dots\times \mathfrak{S}_{\mu_l} \subseteq  \mathfrak{S}_s$.  
As a first step towards understanding the stable Kronecker coefficients,     it is natural to  consider
$$
\dim_{\CC} \Hom_{\mathfrak{S}_s}({\sf M} (\mu), \Delta^0_s(\nu \setminus \lambda)  ) 
$$ 
 and to attempt to construct a basis   in terms of semistandard (Kronecker) tableaux.

 \begin{defn}\label{sstrd}
Let  $(\lambda,\nu,s) \in \mathscr{P}_{r-s} \times \mathscr{P}_{\leq r} \times  \mathbb{N}$ be
a pair of one-row partitions or a triple of maximal depth.  
Let $\mu = (\mu_1, \mu_2, \ldots , \mu_l)\vdash s$ and 
let $\sts, \stt \in \Std^0_s(\nu \setminus \lambda)$.
\begin{enumerate}
\item For $1\leq k <s$ we write $\sts \overset{k}{\sim} \stt$ if $\sts = \stt_{k\leftrightarrow k+1}$.
\item We write $\sts \overset{\mu}{\sim}  \stt$ if there exists a sequence of standard Kronecker tableaux $\stt_1, \stt_2, \ldots , \stt_d \in \Std^0_s(\nu\setminus\lambda)$ 
such that 
$$\sts =  \stt_{1}  \overset{k_1}{\sim}   \stt_{2}  ,  \
 \stt_{2}  \overset{k_2}{\sim}    \stt_{3}  ,   \ \dots \ , 
\stt_{d-1}\overset{k_{d-1}}{\sim}   \stt_{d} 
=\stt  $$  
for some $k_1,\dots, k_{d-1}\in  \{1, \ldots , s-1\} \setminus 
 \{    [\mu]_c \mid  c = 1, \ldots , l-1 \}$.   
We define a {\sf tableau of weight} $\mu$ to be an equivalence class of tableau under $\overset{\mu}{\sim}  $, denoted $[\stt]_\mu = \{ \sts \in \Std^0_s(\nu\setminus \lambda) \, |\, \sts \overset{\mu}{\sim} \stt\}$.
\item We say that a    Kronecker tableau, $[\stt]_\mu$,  of shape $\nu\setminus \lambda$ and weight $\mu$ is {\sf semistandard} if  for any $\sts \in [\stt]_\mu$ and any $  k \not \in  \{[\mu_c] \mid c = 1, \ldots , l-1 \}$ the tableau $\sts_{k\leftrightarrow k+1}$ exists. We 
let $\SStd_s^0(\nu\setminus \lambda, \mu)$ denote the set of  semistandard Kronecker tableaux of shape $\nu\setminus \lambda$ and weight $\mu$.
\end{enumerate}
\end{defn}

\begin{figure}[ht!]
\scalefont{0.75}   
    \begin{tikzpicture}[scale=0.6]
  \draw[white] [decorate,decoration={brace,amplitude=6pt},xshift=6pt,yshift=0pt]   (-3.85,-7.9)--   (-3.85,-0.1)  node [right,black,midway,yshift=-0.2cm,xshift=-2cm]{\text{1st frame}} 
  node [right,black,midway,yshift=0.2cm,xshift=-2cm]{\text{2 steps in}} ; 
  
  \draw[white] [decorate,decoration={brace,amplitude=6pt},xshift=6pt,yshift=0pt]   (-3.85,-7.9-8)--   (-3.85,-0.1-8)  node [right,black,midway,yshift=-0.2cm,xshift=-2cm]{\text{2nd frame}} 
  node [right,black,midway,yshift=0.2cm,xshift=-2cm]{\text{2 steps in}} ; 

  \draw[white] [decorate,decoration={brace,amplitude=6pt},xshift=6pt,yshift=0pt]   (-3.85,-7.9-8-4)--   (-3.85,-0.1-8-8)  node [right,black,midway,yshift=-0.2cm,xshift=-2cm]{\text{3rd frame}} 
  node [right,black,midway,yshift=0.2cm,xshift=-2cm]{\text{1 step in}} ; 

     \clip(-3.6,0.5) rectangle (3.6,-20.8);
    \path   (0,0)edge[decorate]  node[left] {$-1$}  (0,-2);
          \path   (0,-2)edge[decorate]  node[right] {$+1$}  (2,-4);
                    \path   (0,-2)edge[decorate]  node[left] {$+0$}  (-2,-4);

                    \path   (-2,-6)edge[decorate]  node[left] {$-1$}  (-2,-4);                         
          \path   (2,-6)edge[decorate]  node[right] {$-1$}  (2,-4);
                                       
                  \path   (-2,-6)edge[decorate]  node[left] {$+1$}  (0,-8);                         
          \path   (2,-6)edge[decorate]  node[right] {$+0$}  (0,-8);

       \path   (0,0-8)edge[decorate]  node[left] {$-1$}  (-2,-2-8);
              \path   (0,0-8)edge[decorate]  node[right] {$-0$}  (2,-2-8);

          \path   (2,-2-8)edge[decorate]  node[right] {$+1$}  (2,-4-8);
                    \path   (-2,-2-8)edge[decorate]  node[left] {$+1$}  (-2,-4-8);

                    \path   (0,-6-8)edge[decorate]  node[left] {$-0$}  (-2,-4-8);                         
          \path   (0,-6-8)edge[decorate]  node[right] {$-1$}  (2,-4-8);
                                       
                  \path   (0,-6-8)edge[decorate]  node[left] {$+1$}  (0,-8-8);

       \path   (0,-16)edge[decorate]  node[left] {$-0$}  (0,-18);
       \path   (0,-20)edge[decorate]  node[left] {$+1$}  (0,-18);

              \draw[dashed] (-3.5,0) rectangle (3.5,-8);                   
                \draw[dashed] (-3.5,-8)--(-3.5,-16)--(3.5,-16)--(3.5,-8) ;  
                
                      \draw[dashed] (-3.5,-16)--(-3.5,-20)--(3.5,-20)--(3.5,-16) ;  
                  \fill[white] (0,0) circle (17pt);   
  \begin{scope}   
 \fill[white] (0,0.1) circle (22pt);
     \draw (0,0) node {$  \scalefont{0.4}\yng(4) $  };   
\fill[white] (0,-2) circle (17pt);    \draw (0,-2) node{$  \scalefont{0.4}\yng(3)$  };
\fill[white] (-2,-4) circle (17pt);  \draw (-2,-4) node{$  \scalefont{0.4}\yng(3) $  };
\fill[white] (2,-4) circle (17pt);  \draw (2,-4) node{$   \scalefont{0.4}\yng(4) $  };
  \fill[white] (-2,-6) circle (17pt);  \draw (-2,-6) node{$  \scalefont{0.4}\yng(2)$  };
\fill[white] (2,-6) circle (17pt);  \draw (2,-6) node{$  \scalefont{0.4}\yng(3) $  };
  \fill[white] (0,-8) circle (17pt);  \draw (0,-8) node{$  \scalefont{0.4}\yng(3)$  };
  \fill[white] (0,-16) circle (17pt);  \draw (0,-16) node{$  \scalefont{0.4}\yng(4)$  };
\fill[white] (-2,-2-8) circle (17pt);    \draw (-2,-2-8) node{$  \scalefont{0.4}\yng(2)$  };
\fill[white] (2,-2-8) circle (17pt);    \draw (2,-2-8) node{$  \scalefont{0.4}\yng(3)$  };
\fill[white] (-2,-4-8) circle (17pt);  \draw (-2,-4-8) node{$  \scalefont{0.4}\yng(3) $  };
\fill[white] (2,-4-8) circle (17pt);  \draw (2,-4-8) node{$   \scalefont{0.4}\yng(4) $  };
  \fill[white] (0,-6-8) circle (17pt);  \draw (0,-6-8) node{$  \scalefont{0.4}\yng(3)$  };

   \fill[white] (0.5,-20) circle (12pt); \fill[white] (-0.5,-20) circle (12pt);

  \fill[white] (0,-18) circle (17pt);  \draw (0,-18) node{$  \scalefont{0.4}\yng(4)$  };
\fill[white] (0,-20) circle (17pt);  \draw (0,-20) node{$   \scalefont{0.4}\yng(5)$  };

      \end{scope}
   \end{tikzpicture} 
   \qquad 
    \begin{tikzpicture}[scale=0.6]
      \clip(-3.6,0.5) rectangle (3.6,-20.8);
    \path   (0,0)edge[decorate]  node[left] {$-1$}  (-2,-2);
    \path   (0,0)edge[decorate]  node[right] {$-4$}  (2,-2);

          \path   (2,-2)edge[decorate]  node[right] {$+0$}  (2,-4);
                    \path   (-2,-2)edge[decorate]  node[right] {$+0$}  (-2,-4);

                    \path   (-2,-6)edge[decorate]  node[right] {$-4$}  (-2,-4);                         
          \path   (2,-6)edge[decorate]  node[right] {$-1$}  (2,-4);
                                       
                  \path   (-2,-6)edge[decorate]  node[left] {$+0$}  (0,-8);                         
          \path   (2,-6)edge[decorate]  node[right] {$+0$}  (0,-8);

       \path   (0,0-8)edge[decorate]  node[left] {$-1$}  (-2,-2-8);
              \path   (0,0-8)edge[decorate]  node[right] {$-2$}  (2,-2-8);

          \path   (2,-2-8)edge[decorate]  node[right] {$+3$}  (2,-4-8);
                    \path   (-2,-2-8)edge[decorate]  node[left] {$+0$}  (-2,-4-8);

                    \path   (-2,-6-8)edge[decorate]  node[left] {$-2$}  (-2,-4-8);                         
          \path   (2,-6-8)edge[decorate]  node[right] {$-1$}  (2,-4-8);
                                       
                  \path   (-2,-6-8)edge[decorate]  node[left] {$+3$}  (0,-8-8);                         
                  \path   (2,-6-8)edge[decorate]  node[left] {$+0$}  (0,-8-8);

       \path   (0,-16)edge[decorate]  node[right] {$-2$}  (0,-18);
       \path   (0,-20)edge[decorate]  node[right] {$+3$}  (0,-18);

              \draw[dashed] (-3.5,0) rectangle (3.5,-8);                   
                \draw[dashed] (-3.5,-8)--(-3.5,-16)--(3.5,-16)--(3.5,-8) ;  
                
                      \draw[dashed] (-3.5,-16)--(-3.5,-20)--(3.5,-20)--(3.5,-16) ;  
                  \fill[white] (0,0) circle (17pt);   
  \begin{scope}   
 \fill[white] (0,0.1) circle (26pt); \fill[white] (-.2,-0.1) circle (7pt);\fill[white] (1.2,0) circle (12pt);
 \fill[white] (1,-8) circle (12pt); \fill[white] (-1,-8) circle (12pt);
  \fill[white] (1,-16) circle (12pt); \fill[white] (-1,-16) circle (12pt);
   \fill[white] (0.8,-20) circle (12pt); \fill[white] (-0.8,-20) circle (12pt);
 \fill[white] (-0.6,-0.7) circle (7pt);    \draw (0.4,-0.4) node {$  \scalefont{0.4}\yng(7,5,1,1) $  };   
\fill[white] (2,-2) circle (17pt);    \draw (2,-2.2) node{$  \scalefont{0.4}\yng(7,5,1)$  };
\fill[white] (-2,-2) circle (17pt);    \draw (-2,-2.3) node{$  \scalefont{0.4}\yng(6,5,1,1)$  };
\fill[white] (2,-4) circle (17pt);    \draw (2,-4.2) node{$  \scalefont{0.4}\yng(7,5,1)$  };
\fill[white] (-2,-4) circle (17pt);    \draw (-2,-4.3) node{$  \scalefont{0.4}\yng(6,5,1,1)$  };
  \fill[white] (-2,-6.1) circle (18pt);    \draw (-2,-6.1) node{$  \scalefont{0.4}\yng(6,5,1)$  };
\fill[white] (2,-6.1) circle (18pt);    \draw (1.84,-6.1) node{$  \scalefont{0.4}\yng(6,5,1)$  };
  \fill[white] (0,-8) circle (17pt);  
    \fill[white] (0-0.5,-8-0.5) circle (9pt);  \draw (0.1,-8.2) node{$  \scalefont{0.4}\yng(6,5,1)$  };
  \fill[white] (0,-16) circle (17pt);  \draw (0,-16.2) node{$    \scalefont{0.4}\yng(5,4,2)$   };
\fill[white] (-2,-2-8) circle (17pt);    \draw (-2,-2-8.2) node{$   \scalefont{0.4}\yng(5,5,1)$   };
\fill[white] (2,-2-8) circle (17pt);    \draw (2,-2-8.2) node{$  \scalefont{0.4}\yng(6,4,1)$   };
\fill[white] (-2,-4-8) circle (17pt);  \draw (-2,-4-8.1) node{$   \scalefont{0.4}\yng(5,5,1) $  };
\fill[white] (2,-4-8) circle (17pt);  \draw (2,-4-8.1) node{$    \scalefont{0.4}\yng(6,4,2)$   };
\fill[white] (-2,-4-8-2) circle (17pt);  \draw (-2,-4-8.1-2) node{$   \scalefont{0.4}\yng(5,4,1)$    };
\fill[white] (2,-4-8.1-2) circle (18pt);  \draw (2,-4-8.1-2) node{$    \scalefont{0.4}\yng(5,4,2)$   };

  \fill[white] (0,-18) circle (17pt);  \draw (0,-18.2) node{$    \scalefont{0.4}\yng(5,3,2)$    };
\fill[white] (0,-20.2) circle (19pt);  \draw (0,-20.1) node{$    \scalefont{0.4}\yng(5,3,3)$    };

      \end{scope}
   \end{tikzpicture} 

   \caption{Two examples of semistandard Kronecker tableaux of weight $\mu=(2,2,1)$. 
   The number of steps in the $i$th  frame is   $\mu_i$.  
  }
   \label{anewfigforintro}
\end{figure}
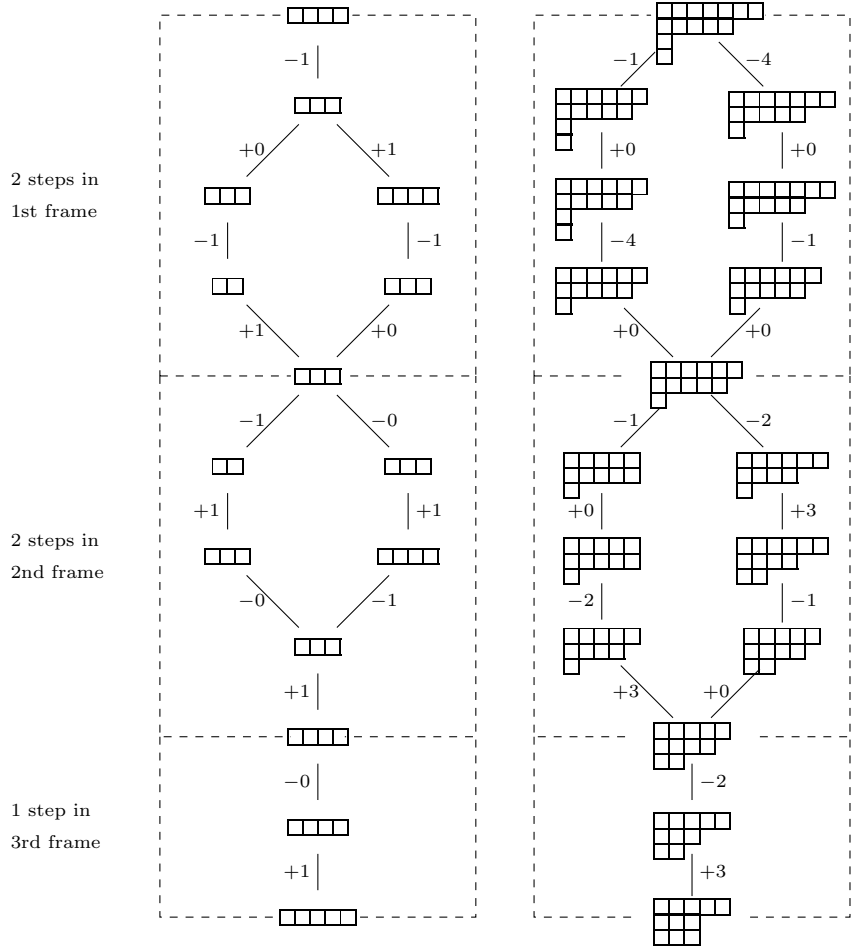

To represent these semistandard Kronecker tableaux graphically, we will add \lq frames' corresponding to the composition $\mu$ on  the set of paths $\Std_s^0(\nu \setminus \lambda)$ in $\mathcal{Y}$. 
For $\stt =(-\varepsilon_{i_1}, + \varepsilon_{j_1}, \ldots , -\varepsilon_{i_s}, + \varepsilon_{j_s})$ we say that the integral step $(-\varepsilon_{i_k}, + \varepsilon_{j_k})$ belongs to the $c$th frame if $[\mu]_{c-1} < k\leq  [\mu]_c$.
Thus for $\sts, \stt \in \Std_s^0(\nu\setminus \lambda)$ we have that $\sts \overset{\mu}{\sim} \stt$ if and only if $\sts$ is obtained from $\stt$ by permuting integral steps within each frame (as in Figures   \ref{anewfigforintro} and \ref{maximaldepth}).

  \begin{thm} \label{YOUNGSRULE}
Let $(\lambda, \nu, s)$ be a co-Pieri triple and $\mu\vdash s$. We define 
$ \varphi_\SSTT({\stt^\mu}) = \sum_{\sts\in  \SSTT}\sts   $ 
for $\SSTT \in \SStd_s^0(\nu \setminus \lambda, \mu)$.   Then 
$   \Hom_{\mathfrak{S}_s}({\sf M}(\mu), \Delta_s^0(\nu\setminus \lambda))$ has 
$\ZZ$-basis $ \{\varphi_\SSTT \mid \SSTT\in  \SStd_s^0(\nu \setminus \lambda, \mu)\}$.    
 \end{thm}

%

 \begin{eg}\label{semiexam2}
  Let  
  $\lambda =(4) $, $\nu =(4)$ and $s=5$ and  $\mu=(2,2,1) \vdash {5}$.   
An example of a semistandard  tableau, $\SSTV $,   of shape $\nu\setminus \lambda$  and weight $\mu$ is given by the rightmost diagram in  Figure \ref{anewfigforintro}.  The semistandard tableau $\SSTV$ is an orbit consisting of the following four standard tableaux
\begin{align*}
&\stv_1= r(1) \circ d(1) \circ d(1) \circ a(1) \circ a(1)
 \ \quad  \ 
\stv_2=   d(1) \circ r(1) \circ d(1) \circ a(1) \circ a(1)
  \  \quad  \ \\
& \stv_3=   r(1) \circ d(1) \circ a(1) \circ d(1) \circ a(1)
 \  \quad  \ 
 \stv_4= d(1) \circ r(1) \circ a(1) \circ d(1) \circ a(1)
 \end{align*}
 We have a corresponding homomorphism
 $ 
 \varphi_\SSTV \in \Hom_{\mathfrak{S}_s}({\sf M}(2,2,1), \Delta_s((4)\setminus (4))
 $ 
 given by 
 $$\varphi_\SSTT(\stt^{(2,2,1)})=\stv_1+\stv_2+\stv_3+\stv_4.  $$
\end{eg}

\subsection{The classical picture for semistandard Young tableaux}
 We now wish to illustrate how our Definition \ref{sstrd}  and   the familiar visualisation  of a  semistandard  Young tableaux coincide for triples of maximal depth.  
%
%
      Given   $\lambda \vdash  {r-s} , \nu \vdash  {r}, \mu = (\mu_1, \mu_2, \ldots , \mu_\ell ) \vdash s $
   such that  $\lambda \subseteq \nu$   a   Young tableau of shape $\nu\ominus \lambda$ and weight $\mu$   in the classical picture is visualised as   a filling of the boxes of 
   $[\nu\ominus \lambda]$  with the entries 
    $$\underbrace{1, \dots, 1}_{\mu_1}, \underbrace{2,\dots, 2}_{\mu_2}, 
  \ldots,   \underbrace{\ell ,\dots, \ell }_{\mu_\ell } $$ 
 so that they are weakly increasing along the rows and columns.  
    One should think of this classical picture of a 
       Young tableau of weight $\mu$  
  simply  as a  diagrammatic way of encoding an
  $\mathfrak{S}_\mu$-orbit  of   standard  Young tableaux as follows.  
  Let $\sts$ be a standard   Young tableau of shape $\nu\ominus \lambda$    and let $\mu$  be a partition. Then define $\mu(\sts)$ to be  the  Young tableau of weight $\mu$
 obtained from $\sts$ by replacing each of the entries
  $[\mu]_{c-1} < i \leq [\mu]_c$ in $\sts$ by the entry $c$  for  $  c \geq  1$.  
We identify a   Young tableau, $\SSTS $, of weight $\mu$ with the set of standard  Young tableaux, $\mu^{-1}(\SSTS)=\{\sts \mid \mu(\sts)=\SSTS\}$. 
  
 In either picture, a 
   Young tableau of weight $\mu$
  is merely a picture which  encodes an $\mathfrak{S}_\mu$-orbit of standard  Young tableaux. 
     We  
  picture a   Young tableau, $\SSTS$, of weight $\mu$   as the   orbit of  paths $\mu^{-1}(\SSTS)$ in the branching graph with a frame to record the partition $\mu$. 

   A  tableau of weight $\mu$ in the classical picture
    would be said to be semistandard if and only if the entries are strictly increasing along the columns.  
In our picture, this is   equivalent to condition 3 of Definition \ref{sstrd}.

  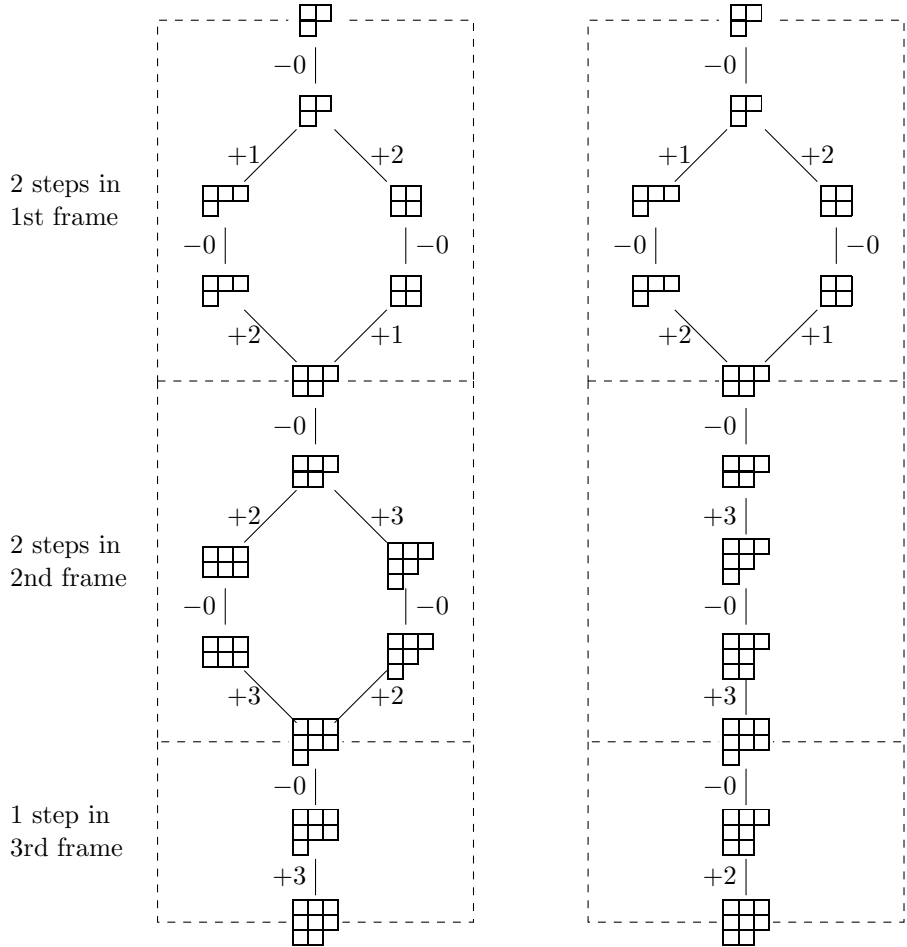
\begin{figure}[ht!]
 $$ \begin{tikzpicture}[scale=0.6]
  \draw[white] [decorate,decoration={brace,amplitude=6pt},xshift=6pt,yshift=0pt]   (-3.85,-7.9)--   (-3.85,-0.1)  node [right,black,midway,yshift=-0.2cm,xshift=-2cm]{\text{1st frame}} 
  node [right,black,midway,yshift=0.2cm,xshift=-2cm]{\text{2 steps in}} ; 
  
  \draw[white] [decorate,decoration={brace,amplitude=6pt},xshift=6pt,yshift=0pt]   (-3.85,-7.9-8)--   (-3.85,-0.1-8)  node [right,black,midway,yshift=-0.2cm,xshift=-2cm]{\text{2nd frame}} 
  node [right,black,midway,yshift=0.2cm,xshift=-2cm]{\text{2 steps in}} ; 

  \draw[white] [decorate,decoration={brace,amplitude=6pt},xshift=6pt,yshift=0pt]   (-3.85,-7.9-8-4)--   (-3.85,-0.1-8-8)  node [right,black,midway,yshift=-0.2cm,xshift=-2cm]{\text{3rd frame}} 
  node [right,black,midway,yshift=0.2cm,xshift=-2cm]{\text{1 step in}} ; 

  \clip(-3.6,0.5) rectangle (3.6,-20.8);
    \path   (0,0)edge[decorate]  node[left] {$-0$}  (0,-2);
          \path   (0,-2)edge[decorate]  node[right] {$+2$}  (2,-4);
                    \path   (0,-2)edge[decorate]  node[left] {$+1$}  (-2,-4);

                    \path   (-2,-6)edge[decorate]  node[left] {$-0$}  (-2,-4);                         
          \path   (2,-6)edge[decorate]  node[right] {$-0$}  (2,-4);
                                       
                  \path   (-2,-6)edge[decorate]  node[left] {$+2$}  (0,-8);                         
          \path   (2,-6)edge[decorate]  node[right] {$+1$}  (0,-8);

       \path   (0,0-8)edge[decorate]  node[left] {$-0$}  (0,-2-8);
          \path   (0,-2-8)edge[decorate]  node[right] {$+3$}  (2,-4-8);
                    \path   (0,-2-8)edge[decorate]  node[left] {$+2$}  (-2,-4-8);

                    \path   (-2,-6-8)edge[decorate]  node[left] {$-0$}  (-2,-4-8);                         
          \path   (2,-6-8)edge[decorate]  node[right] {$-0$}  (2,-4-8);
                                       
                  \path   (-2,-6-8)edge[decorate]  node[left] {$+3$}  (0,-8-8);                         
          \path   (2,-6-8)edge[decorate]  node[right] {$+2$}  (0,-8-8);

       \path   (0,-16)edge[decorate]  node[left] {$-0$}  (0,-18);
       \path   (0,-20)edge[decorate]  node[left] {$+3$}  (0,-18);

              \draw[dashed] (-3.5,0) rectangle (3.5,-8);                   
                \draw[dashed] (-3.5,-8)--(-3.5,-16)--(3.5,-16)--(3.5,-8) ;  
                
                      \draw[dashed] (-3.5,-16)--(-3.5,-20)--(3.5,-20)--(3.5,-16) ;  
                  \fill[white] (0,0) circle (17pt);   
  \begin{scope}   
 \fill[white] (0,0) circle (17pt);
     \draw (0,0) node {$  \scalefont{0.4}\yng(2,1) $  };   
\fill[white] (0,-2) circle (17pt);    \draw (0,-2) node{$  \scalefont{0.4}\yng(2,1)$  };
\fill[white] (-2,-4) circle (17pt);  \draw (-2,-4) node{$  \scalefont{0.4}\yng(3,1) $  };
\fill[white] (2,-4) circle (17pt);  \draw (2,-4) node{$   \scalefont{0.4}\yng(2,2) $  };
  \fill[white] (-2,-6) circle (17pt);  \draw (-2,-6) node{$  \scalefont{0.4}\yng(3,1)$  };
\fill[white] (2,-6) circle (17pt);  \draw (2,-6) node{$  \scalefont{0.4}\yng(2,2) $  };
  \fill[white] (0,-8) circle (17pt);  \draw (0,-8) node{$  \scalefont{0.4}\yng(3,2)$  };
  \fill[white] (0,-16) circle (17pt);  \draw (0,-16) node{$  \scalefont{0.4}\yng(3,3,1)$  };
\fill[white] (0,-2-8) circle (17pt);    \draw (0,-2-8) node{$  \scalefont{0.4}\yng(3,2)$  };
\fill[white] (-2,-4-8) circle (17pt);  \draw (-2,-4-8) node{$  \scalefont{0.4}\yng(3,3) $  };
\fill[white] (2,-4-8) circle (17pt);  \draw (2.1,-4-8.1) node{$   \scalefont{0.4}\yng(3,2,1) $  };
  \fill[white] (-2,-6-8) circle (17pt);  \draw (-2,-6-8) node{$  \scalefont{0.4}\yng(3,3)$  };
\fill[white] (2,-6-8) circle (17pt);  \draw (2.1,-6-8.1) node{$   \scalefont{0.4}\yng(3,2,1)$  };

  \fill[white] (0,-18) circle (17pt);  \draw (0,-18) node{$  \scalefont{0.4}\yng(3,3,1)$  };
\fill[white] (0,-20) circle (17pt);  \draw (0,-20) node{$   \scalefont{0.4}\yng(3,3,2)$  };

      \end{scope}
   \end{tikzpicture} 
   \qquad\qquad
   \begin{tikzpicture}[scale=0.6]
   
  \clip(-3.6,0.5) rectangle (3.6,-20.8);
    \path   (0,0)edge[decorate]  node[left] {$-0$}  (0,-2);
          \path   (0,-2)edge[decorate]  node[right] {$+2$}  (2,-4);
                    \path   (0,-2)edge[decorate]  node[left] {$+1$}  (-2,-4);

                    \path   (-2,-6)edge[decorate]  node[left] {$-0$}  (-2,-4);                         
          \path   (2,-6)edge[decorate]  node[right] {$-0$}  (2,-4);
                                       
                  \path   (-2,-6)edge[decorate]  node[left] {$+2$}  (0,-8);                         
          \path   (2,-6)edge[decorate]  node[right] {$+1$}  (0,-8);

       \path   (0,0-8)edge[decorate]  node[left] {$-0$}  (0,-2-8);
          \path   (0,-2-8)edge[decorate]  node[left] {$+3$}  (0,-4-8);

                    \path   (0,-6-8)edge[decorate]  node[left] {$-0$}  (0,-4-8);                         
                                       
                  \path   (0,-6-8)edge[decorate]  node[left] {$+3$}  (0,-8-8);

       \path   (0,-16)edge[decorate]  node[left] {$-0$}  (0,-18);
       \path   (0,-20)edge[decorate]  node[left] {$+2$}  (0,-18);

              \draw[dashed] (-3.5,0) rectangle (3.5,-8);                   
                \draw[dashed] (-3.5,-8)--(-3.5,-16)--(3.5,-16)--(3.5,-8) ;  
                
                      \draw[dashed] (-3.5,-16)--(-3.5,-20)--(3.5,-20)--(3.5,-16) ;  
                  \fill[white] (0,0) circle (17pt);   
  \begin{scope}   
 \fill[white] (0,0) circle (17pt);
     \draw (0,0) node {$  \scalefont{0.4}\yng(2,1) $  };   
\fill[white] (0,-2) circle (17pt);    \draw (0,-2) node{$  \scalefont{0.4}\yng(2,1)$  };
\fill[white] (-2,-4) circle (17pt);  \draw (-2,-4) node{$  \scalefont{0.4}\yng(3,1) $  };
\fill[white] (2,-4) circle (17pt);  \draw (2,-4) node{$   \scalefont{0.4}\yng(2,2) $  };
  \fill[white] (-2,-6) circle (17pt);  \draw (-2,-6) node{$  \scalefont{0.4}\yng(3,1)$  };
\fill[white] (2,-6) circle (17pt);  \draw (2,-6) node{$  \scalefont{0.4}\yng(2,2) $  };
  \fill[white] (0,-8) circle (17pt);  \draw (0,-8) node{$  \scalefont{0.4}\yng(3,2)$  };
  \fill[white] (0,-16) circle (17pt);  \draw (0,-16) node{$  \scalefont{0.4}\yng(3,3,1)$  };
\fill[white] (0,-2-8) circle (17pt);    \draw (0,-2-8) node{$  \scalefont{0.4}\yng(3,2)$  };
\fill[white] (0,-4-8) circle (17pt);  \draw (0,-4-8) node{$  \scalefont{0.4}\yng(3,2,1) $  };
\fill[white] (0,-6-8) circle (17pt);  \draw (0,-6-8.1) node{$   \scalefont{0.4}\yng(3,2,2)$  };

  \fill[white] (0,-18) circle (17pt);  \draw (0,-18) node{$  \scalefont{0.4}\yng(3,2,2)$  };
\fill[white] (0,-20) circle (17pt);  \draw (0,-20) node{$   \scalefont{0.4}\yng(3,3,2)$  };

      \end{scope}
   \end{tikzpicture} $$
   \caption{A pair of semistandard Kronecker tableaux for $((2,1),(3,3,2),5)$ a triple of maximal depth.  Compare the first of these with \cref{semiexam1}.}
   \label{maximaldepth}
   \end{figure}

 \begin{eg}\label{semiexam1}
  Let  
  $\lambda =(2,1) $, $\nu =(3,3,2)$ and $s=5$. Then $(\lambda, \nu,s)$ is a triple of maximal depth. Take $\mu=(2,2,1) \vdash {5}$.   
 The semistandard tableau $\SSTU$ is an orbit consisting of the following four standard tableaux
\begin{align*}
&\stu_1= a(1) \circ a(2) \circ a(2) \circ a(3) \circ a(3)
 \ \quad  \ 
\stu_2=  a(2) \circ a(1) \circ a(2) \circ a(3) \circ a(3)
  \  \quad  \ \\
& \stu_3=  a(1) \circ a(2) \circ a(3) \circ a(2) \circ a(3)
 \  \quad  \ 
 \stu_4= a(2) \circ a(1) \circ a(3) \circ a(2) \circ a(3) 
 \end{align*}
 pictured as follows
 $$\scalefont{0.9}
\Yboxdim{12pt}
\mu^{-1} \left( \;  \Yvcentermath1 \young(\ \ 1,\ 12,23)\ \right) = 
\left\{
\;  \Yboxdim{12pt} \Yvcentermath1 \young(\ \ 1,\ 23,45)  
\ \ , \ \ 
\Yboxdim{12pt} \Yvcentermath1 \young(\ \ 2,\ 13,45) 
\ \ , \ \ 
\Yboxdim{12pt} \Yvcentermath1 \young(\ \ 2,\ 14,35) 
\ \ , \ \ 
\Yboxdim{12pt} \Yvcentermath1 \young(\ \ 1,\ 23,45) 
\ \right\}.
 $$
    We have a corresponding homomorphism
 $ 
 \varphi_\SSTU \in \Hom_{\mathfrak{S}_s}({\sf M}(2,2,1), \Delta_s((3,3,2)\setminus (2,1))
 $ 
 given by 
 $$\varphi_\SSTT(\stt^{(2,2,1)})=\stu_1+\stu_2+\stu_3+\stu_4.  $$
Compare this orbit sum over 4 tableaux with the   picture in Figure \ref{maximaldepth} and the statement of Theorem \ref{YOUNGSRULE}.  
\end{eg}

\begin{eg}
 Let  
  $\lambda =(2,1) $, $\nu =(3,3,2)$ and $s=5$. Then $(\lambda, \nu,s)$ is a triple of maximal depth. Take $\mu=(2,2,1) \vdash {5}$.   
 The full list of  semistandard tableaux (pictured in the  classical fashion) are as follows 
 $$
 \Yvcentermath1 \young(\ \ 1,\ 12,23)
 \quad
 \Yvcentermath1 \young(\ \ 1,\ 13,22)
 \quad
  \Yvcentermath1 \young(\ \ 1,\ 22,13)
  \quad
   \Yvcentermath1 \young(\ \ 2,\ 13,12)
 $$
 The first two of these semistandard tableaux  are pictured in our diagrammatic fashion in \cref{maximaldepth}.  
\end{eg}

  \section{Latticed Kronecker tableaux}\label{sec:latticed}

We now provide the main result of the paper, namely we  combinatorially describe
$$\overline{g}(\lambda, \nu, \mu) = \dim \Hom_{\mathfrak{S}_s}({\sf S}(\mu), \Delta_s^0(\nu\setminus \lambda))$$
for  $(\lambda, \nu, \mu)$  a triple of maximal depth or such that $\lambda$ and $\nu$ are both one-row partitions.     
 One can think of a   path $\stt \in \Std_s(\nu\setminus\lambda)$ as a sequence of partitions; or equivalently, as the 
  sequence of  
    boxes added 
 and removed.   We shall   refer to a  pair of steps, $(-\varepsilon_a,+\varepsilon_b)$,   between consecutive  integral levels of the branching graph    as an {\sf integral step} in the branching graph.
  We define {\sf  types} of integral step  (move-up, dummy, move-down) 
 in the branching graph of $P_r(n)$ and order them as follows, 
 $$\begin{array}{ccccccccc}
 	 &\text{move-up }			&  		 &\text{dummy }		&  		 &\text{move-down }	
&  		 
									\\
									
  		& (-\varepsilon_p, + \varepsilon_q)&< 		& (-\varepsilon_t, + \varepsilon_t)
&< 		&(-\varepsilon_u, + \varepsilon_v)  		 
\end{array}$$
for $p>q$ and $u< v$;   we refine this to a total  order as follows,
\begin{enumerate}
\item[$({\up })$]  we order $(-\varepsilon_p, + \varepsilon_q)< (-\varepsilon_{p'}, + \varepsilon_{q'}) $ if $q<q'$ or $q=q'$ and $p>p'$; 
 \item[$(d)$]  we order $(-\varepsilon_t, + \varepsilon_t) < (-\varepsilon_{t'}, + \varepsilon_{t'}) $ if $t>t'$;
\item[$({\down })$]  we order $(-\varepsilon_u, + \varepsilon_v)< (-\varepsilon_{u'}, + \varepsilon_{v'})$  if $u>u'$ or $u=u'$ and $v<v'$. 
\end{enumerate}
    We sometimes let $a(i):=\down(0,i)$ (respectively $r(i):=\up(i,0)$) and think of this as  {\sf adding} (respectively  {\sf removing}) a box.  
  We start with any standard tableau $\sts \in \Std_s^0(\nu \setminus \lambda)$ and any $\mu = (\mu_1, \mu_2, \ldots , \mu_l)\vdash s$. Write 
$$\sts =  (-\varepsilon_{i_1}, 
+\varepsilon_{j_1},
-\varepsilon_{i_2},
+\varepsilon_{j_2},
\dots
, -\varepsilon_{i_s},
+\varepsilon_{j_s}). $$
Recall from the previous section that, to each integral step $(-\varepsilon_{i_k}, + \varepsilon_{j_k})$ in $\sts$, we associate its frame  $c$, that is the unique positive integer such that 
 $[\mu]_{c-1} < k \leq [\mu]_c.$

\begin{defn}\label{jdfhklssdhjhlashlfs}   
We encode the integral steps of $\sts$ and their frames in a $2\times s$ array, denoted by $\omega_\mu (\sts)$ (called the $\mu$-reverse reading word of $\sts$) as follows. The first row of $\omega_\mu(\sts)$ contains all the integral steps of $\sts$ and the second row contains their corresponding frames. We order the columns of $\omega_\mu(\sts)$ increasingly using the ordering on integral steps given in Definition 2.5. 
For two equal integral steps we order the columns so that the frame numbers are weakly decreasing.  
 Given $\SSTS\in \SStd_s^0(\nu\setminus \lambda, \mu)$,  it is easy to see that 
 $\omega_\mu (\sts)=\omega_\mu (\stt)$ for any pair $\sts,\stt \in \SSTS$ and so we define the $\mu$-reverse reading word, $\omega(\SSTS)$, of $\SSTS$ in the obvious fashion.    
For $\SSTS\in \SStd_s^0(\nu\setminus \lambda, \mu)$ we write 
$$\omega(\SSTS) = (\omega_1(\SSTS), \omega_2(\SSTS))$$
where $\omega_1(\SSTS)$ (respectively  $\omega_2(\SSTS)$) is the first (respectively  second) row of $\omega(\SSTS)$.
Note that $\omega_2(\SSTS)$ is a sequence   of positive integers such that $i$ appears precisely $\mu_i$ times, for  $i\geq 1$.
\end{defn}

\begin{eg}\label{semiexam3}
For $\lambda=(2,1)$ and $\nu=(3,3,2)$, the steps taken in the semistandard tableau $\SSTU$  of Figure \ref{maximaldepth} are 
$$a(1), a(2), a(2), a(3), a(3)$$ 
  We record the  steps according to the dominance ordering for the partition algebra  
 ($a(1)< a(2) < a(3)$)  and   refine this by recording the  frame   in which these steps  occur backwards,  
  as follows 
$$\omega(\SSTU)=
\left(\begin{array}{cccccccccccc}
a(1)&a(2)&a(2)&a(3) &a(3) 
\\
1& 2& 1 & 3 &2
\end{array}\right).$$
%
  For $\lambda=(4)$ and $\nu=(5)$, the steps taken in the semistandard tableau
 $\SSTV$ on the right of Figure \ref{anewfigforintro} are 
$$r(1), d(1), d(1), a(1), a(1).$$
  We record the  steps according to the  dominance ordering for the partition algebra 
 ($r(1)< d(1) < a(1)$)  and we  refine  this by recording the  frame   in which these steps  occur backwards,  
  as follows 
$$
\omega(\SSTV)=\left(\begin{array}{cccccccccccc}
r(1)&d(1)&d(1)&a(1) &a(1) 
\\
1& 2& 1 & 3 &2
\end{array}\right)$$
 and notice that $\omega_2(\SSTU)=\omega_2(\SSTV)$.  We leave it as an exercise for the reader to verify that the rightmost tableau depicted in Figure \ref{anewfigforintro} has reading word
 $$
 \left(\begin{array}{cccccccccccc}
 \ r(4) \ & \ r(1) \ & \ r(1) \ &m{\downarrow}(2,3) &m{\downarrow}(2,3) 
\\
1& 2& 1 & 2&3
\end{array}\right).$$
 \end{eg}


\begin{thm}
For $\SSTS\in \SStd_s^0(\nu\setminus \lambda, \mu)$ we say that its reverse reading word $\omega(\SSTS)$ is a lattice permutation if $\omega_2(\SSTS)$  
is a string composed of positive integers, in which every prefix contains at least as many positive integers $i$ as integers $i+1$ for $i\geq 1$.  
 We define $\Latt_s^0(\nu \setminus \lambda, \mu)$ to be the set of all $\SSTS\in \SStd_s^0(\nu\setminus \lambda, \mu)$ such that $\omega(\SSTS)$ is a lattice permutation.
%
%
%
%
 For any co-Pieri triple $(\lambda, \nu, s)$ and any $\mu\vdash s$ we have 
$$\overline{g}(\lambda, \nu, \mu) = \dim_{\CC}\Hom_{\mathfrak{S}_s}({\sf S}(\mu), \Delta_s^0(\nu\setminus \lambda)) = |\Latt_s^0(\nu\setminus \lambda, \mu)|.$$
\end{thm}

  \begin{eg}
For example, we have that  
  $$
  \overline{g}((2,1), (3,3,2), (2,2,1))= 1  = 
  \overline{g}((4),(4),(2,2,1))
  $$
  and that the corresponding homomorphisms are constructed in Examples \ref{semiexam2} and \ref{semiexam1}.  
  That these semistandard tableaux satisfy the lattice permutation property is checked in Example \ref{semiexam3}.  
  Verifying that these are the only semistandard tableaux satisfying the lattice permutation property is left as an exercise for the reader.   
    Similarly, one can check that $  \overline{g}((7,5,1^2), (6,3,3), (2,2,1))= 1  $.  
%

  \end{eg}

   \begin{rmk}
The (non-stable) Kronecker coefficients
 are
 also indexed by partitions.  
   As we increase the size of the first row of each of the indexing partitions of the   Kronecker coefficients, 
 we obtain a weakly increasing sequence of coefficients; the   limiting values  of these sequences are the stable Kronecker coefficients   which  have been the focus of this paper.  
The non-stable Kronecker coefficients labelled by two 2-line partitions 
 can be written as  an alternating sum of 
 at most 4 stable Kronecker coefficients    labelled by two 1-line partitions  \cite[Proposition 7.6]{BDE}.   
 (In fact, any non-stable Kronecker coefficient can be written as an alternating sum of stable Kronecker coefficients.)
This should be compared with the existing descriptions of     Kronecker coefficients  labelled by two 2-line partitions  \cite{RW94,Rosas01} which also involve 
alternating sums with at most 4 terms.  

The advantages  of our description are that $(1)$ ours is the first description that generalises to other stable Kronecker coefficients (and in particular the first description of any family of  Kronecker coefficients  subsuming the 
 Littlewood--Richardson coefficients)
 and $(2)$ it counts explicit homomorphisms  and therefore works on a higher structural level than all  other descriptions  of stable Kronecker coefficients since those first considered by Littlewood and Richardson \cite{LR34}.
   \end{rmk}

\bibliographystyle{amsplain}   \bibliography{bib}

 \end{document}